\documentclass{ws-ijmpb}
\usepackage{amssymb}

\def\draftnote{}  

\newcommand{\be}{\begin{equation}}
\newcommand{\ee}{\end{equation}}
\newcommand{\ba}{\begin{eqnarray}}
\newcommand{\ea}{\end{eqnarray}}
\newcommand{\ban}{\begin{eqnarray*}}
\newcommand{\ean}{\end{eqnarray*}}
\newcommand{\n}{\nonumber \\}

\newcommand{\eq}[1]{(\ref{#1})}
\newcommand{\sfrac}[2]{{\textstyle \frac{#1}{#2}}}
\newcommand{\ds}{\displaystyle}
\newcommand{\no}{{\textstyle{\circ\atop\circ}}}
\newcommand{\ket}[1]{{| #1 \rangle}}
\newcommand{\qed}{\hfill \fbox{}\medskip}
\newcommand{\ignore}[1]{}
\newcommand{\Z}{\mathbb{Z}}
\newcommand{\DWA}{{\rm DWA}(\widehat{\mathfrak{sl}}_N\!)}
\newcommand{\ZA}{{\rm ZA}(\widehat{\mathfrak{sl}}_N\!)_k}
\newcommand{\DZA}{{\rm DZA}(\widehat{\mathfrak{sl}}_N\!)_k}

\newcommand{\reftitle}[1]{}
\newcommand{\refj}[1]{{\it #1}}
\newcommand{\hep}[1]{}
\newcommand{\hepp}[1]{#1.}

\begin{document}

\markboth{Satoru Odake}
{Comments on the Deformed $W_N$ Algebra}

\catchline{}{}{}

\title{Comments on the Deformed $W_N$ Algebra\footnote{
Talk at the APCTP-Nankai Joint Symposium on ``Lattice Statistics and 
Mathematical Physics'', 8-10 October 2001, Tianjin China.}
}

\author{\footnotesize SATORU ODAKE\footnote{
odake@azusa.shinshu-u.ac.jp}}

\address{Department of Physics, Faculty of Science, Shinshu University\\
Matsumoto 390-8621, Japan}

\maketitle

\pub{DPSU-01-2}{math.QA/0111230}  

\begin{abstract}
We obtain an explicit expression for the defining relation of the 
deformed $W_N$ algebra, $\DWA_{q,t}$. Using this expression we can 
show that, in the $q\rightarrow 1$ limit, 
$\DWA_{q,t}$ with $t=e^{-\frac{2\pi i}{N}}q^{\frac{k+N}{N}}$ reduces to
the $\mathfrak{sl}_N$-version of the Lepowsky-Wilson's $\cal{Z}$-algebra 
of level $k$, $\ZA$.
In other words $\DWA_{q,t}$ with $t=e^{-\frac{2\pi i}{N}}q^{\frac{k+N}{N}}$ 
can be considered as a $q$-deformation of $\ZA$.

In the appendix given by H.~Awata, S.~Odake and J.~Shiraishi, 
we present an interesting relation between $\DWA_{q,t}$ and 
$\zeta$-function regularization.
\end{abstract}

\section{Introduction}

One of our motivation for study of elliptic algebras (deformed Virasoro and
$W$ algebras, elliptic quantum groups, etc.) is to clarify the 
symmetry of massive integrable models. Massive integrable models includes
quantum field theories with mass scale and solvable statistical lattice models.
Typical examples of the latter are models based on $\mathfrak{sl}_2$:
Andrews-Baxter-Forester (ABF) model and Baxter's eight vertex model.
About these models we know the following:\cite{O99}

\medskip
\noindent
\newcommand{\lw}[1]{\smash{\lower2.ex\hbox{#1}}}
\begin{tabular}{|c|c|c|}
\hline
\lw{model}&\lw{{\large ABF}(\uppercase\expandafter{\romannumeral3})}&
\lw{{\large 8 vertex}}\\ 
&&\\
\hline
Boltzmann weight&face type&vertex type\\
\hline
\lw{algebra}&${\cal B}_{q,\lambda}(\widehat{\mathfrak{sl}}_2)$&
\lw{${\cal A}_{q,p}(\widehat{\mathfrak{sl}}_2)$}\\
&$\bigl({\cal B}\otimes\{P,e^Q\}=U_{q,p}(\widehat{\mathfrak{sl}}_2)\bigr)$&\\
\hline
gradation&\lw{homogeneous gradation}&\lw{principal gradation}\\
(energy level of $H_C$)&&\\
\hline
space of states&irr. rep. space of DVA&irr. rep. space of 
${\cal A}_{q,p}(\widehat{\mathfrak{sl}}_2)$\\
\hline
free field&direct&indirect\\
realization&(construction of VO)&(map to ABF)\\
\hline
\end{tabular}

In order to obtain more direct free field realization of the eight vertex model
and its higher rank generalization, it may be useful to study (deformed) 
current algebras of $\mathfrak{sl}_N$ in principal gradation.
Motivated by this, 
Hara {\it et.al} \cite{HJKOS00} studied free field realization 
of the Lepowsky-Wilson's ${\cal Z}$-algebra\cite{LW81} and found some 
relation between the deformed Virasoro algebra (DVA) and 
${\cal Z}$-algebra.
Recently Shiraishi constructed a direct free field realization of the eight 
vertex model with a specific parameter $p=q^3$, 
where the type \uppercase\expandafter{\romannumeral2} vertex operator is 
given by the DVA current.\cite{S01}

In this article we extend the relation between DVA and ${\cal Z}$-algebra
to the higher rank case. In section 2 we present an explicit expression 
for the defining relation of the deformed $W_N$ algebra. 
This is a main result of this article.
In section 3, by using this explicit expression, we show that 
the deformed $W_N$ algebra reduces to the $\mathfrak{sl}_N$-version of the 
Lepowsky-Wilson's $\cal{Z}$-algebra in some limit.
In the appendix given by Awata, Odake and Shiraishi, 
we present an interesting relation between the deformed $W_N$ algebra and 
$\zeta$-function regularization.\cite{AOS00}

\section{Deformed $W_N$ Algebra}

\subsection{Definition}
Let us recall the definition of the deformed $W_N$ algebra, $\DWA_{q,t}$
.\cite{FF95,AKOS95}
It is defined through a free field realization.
This algebra has two parameters($q$ and $t$), and we set $t=q^{\beta}$ and 
$p=qt^{-1}$.
Let us introduce fundamental bosons $h^i_n$ ($n\in\Z$ ; $i=1,\cdots,N$ ; 
$\sum_{i=1}^Np^{in}h^i_n=0$) which satisfy
\be
  [h^i_n,h^j_m]=-\frac{1}{n}(1-q^n)(1-t^{-n})
  \frac{1-p^{(N\delta_{i,j}-1)n}}{1-p^{Nn}}p^{Nn\theta(i<j)}\delta_{n+m,0},
\ee
where $\theta(P)=1$ or $0$ if the proposition $P$ is true or false, 
respectively.
Exponentiated boson $\Lambda_i(z)$ ($i=1,\cdots,N$) is defined by
\be
  \Lambda_i(z)=\;:\exp\Bigl(\sum_{n\neq 0}h^i_nz^{-n}\Bigr)\!:
  q^{\sqrt{\beta}h^i_0}p^{\frac{N+1}{2}-i}.
\ee
Here $:*:$ stands for the usual normal ordering for bosons,
{\it i.e.}, $h^i_n$ with $n\geq 0$ are in the right.
By using this $\Lambda_i(z)$, $\DWA_{q,t}$ current 
$W^i(z)=\sum_{n\in\Z}W^i_nz^{-n}$ ($i=1,\cdots,N-1$) is given by
\be
  W^i(z)=\sum_{1\leq j_1<j_2<\cdots<j_i\leq N}
  :\!\Lambda_{j_1}(p^{\frac{i-1}{2}}z)
   \Lambda_{j_2}(p^{\frac{i-3}{2}}z)\cdots
   \Lambda_{j_i}(p^{-\frac{i-1}{2}}z)\!:\,,
\ee
and we set $W^0(z)=W^N(z)=1$.
(Remark that $\Lambda_i(z)$ corresponds to the weight of vector representation 
of $\mathfrak{sl}_N$ and $W^i(z)$ corresponds to the $i$-th rank antisymmetric
tensor representation.)
$\DWA_{q,t}$ is defined as an associative algebra over $\mathbb{C}$ 
generated by $W^i_n$.\footnote{
It is also defined as a commutant of the screening currents.\cite{FF95,FR97}}

The highest weight state $\ket{\lambda}$ is characterized by 
$W^i_n\ket{\lambda}=0$ ($n>0$) and $W^i_0\ket{\lambda}=
w^i(\lambda)\ket{\lambda}$ ($w^i(\lambda)\in\mathbb{C}$), and the highest 
weight representation space is obtained by successive action of $W^i_{-n}$ 
($n>0$).

Since $\DWA_{q,t}$ has two parameters($q$ and $t$), we can take its
various limit by relating $q$ and $t$.
In the following limit\footnote{
Usually we call this limit as a conformal limit. 
However there are many other limits in which the resultant algebras are
related to conformal field theory.}
\be
  \mbox{Limit I : }
  \biggl\{\begin{array}{ll}
  q=e^{\hbar}\,,\quad&\hbar\rightarrow 0\\
  t=q^{\beta}\,,&\beta\mbox{ : fixed}\quad
  (\alpha_0=\sqrt{\beta}-\frac{1}{\sqrt{\beta}})\, ,
  \end{array}
  \label{LimitI}
\ee
$\DWA_{q,t}$ reduces to the $W_N$ algebra with the Virasoro central charge 
$c=(N-1)(1-N(N+1)\alpha_0^2)$ 
because the $q$-Miura transformation of $\DWA$ becomes the Miura 
transformation of $W_N$ algebra.
Each DWA current $W^i(z)$, however, reduces to some linear combination of
$W_N$ currents.

\subsection{Relation}
In order to write down relations between DWA currents, we define 
the delta function $\delta(z)=\sum_{n\in\Z}z^n$ and the structure
function $f^{i,j}(z)=\sum_{\ell\geq 0}f^{i,j}_{\ell}z^{\ell}$ 
($1\leq i,j\leq N-1$),
\be
  f^{i,j}(z)=\exp\biggl(\sum_{n>0}\frac{1}{n}(1-q^n)(1-t^{-n})
  \frac{1-p^{\min(i,j)n}}{1-p^n}\frac{1-p^{(N-\max(i,j))n}}{1-p^{Nn}}
  p^{\frac{|i-j|}{2}n}z^n\biggr).
  \label{fij}
\ee
It has been expected that DWA currents satisfy quadratic relations, 
$f^{i,j}(\frac{z_2}{z_1})$$W^i(z_1)$ $W^j(z_2)$
$-W^j(z_2)W^i(z_1)f^{j,i}(\frac{z_1}{z_2})=(\mbox{terms containing
delta function})$, in mode expansion 
it becomes
\ba
  [W^i_n,W^j_m]&=&-\sum_{\ell\geq 1}f^{i,j}_{\ell}\bigl(
  W^i_{n-\ell}W^j_{m+\ell}-W^j_{m-\ell}W^i_{n+\ell}\bigr)\n
  &&+(\mbox{contribution from the terms containing delta function}).
  \label{fWW}
\ea
For $i=1$ and $j\geq 1$ case, the relation is\cite{FF95,AKOS95}
\ba
  \!\!\!\!\!\!\!\!\!\!&&f^{1,j}(\sfrac{z_2}{z_1})W^1(z_1)W^j(z_2)
  -W^j(z_2)W^1(z_1)f^{j,1}(\sfrac{z_1}{z_2})\qquad(j\geq 1)\n
  \!\!\!\!\!\!\!\!\!\!&=&
  -\frac{(1-q)(1-t^{-1})}{1-p}\Bigl(
  \delta(p^{\frac{j+1}{2}}\sfrac{z_2}{z_1})W^{j+1}(p^{\frac12}z_2)
  -\delta(p^{-\frac{j+1}{2}}\sfrac{z_2}{z_1})W^{j+1}(p^{-\frac12}z_2)
  \Bigr),
  \label{W1Wj}
\ea
and for $i=2$ and $j\geq 2$ case, the relation is\cite{AKOS95}
\ba
  \!\!\!\!\!\!\!\!\!\!&&f^{2,j}(\sfrac{z_2}{z_1})W^2(z_1)W^j(z_2)
  -W^j(z_2)W^2(z_1)f^{j,2}(\sfrac{z_1}{z_2})\qquad(j\geq 2)\n
  \!\!\!\!\!\!\!\!\!\!&=&
  -\frac{(1-q)(1-t^{-1})}{1-p}\frac{(1-qp)(1-t^{-1}p)}{(1-p)(1-p^2)}\n
  \!\!\!\!\!\!\!\!\!\!&&\qquad\qquad\times
  \Bigl(\delta(p^{\frac{j}{2}+1}\sfrac{z_2}{z_1})W^{j+2}(pz_2)
  -\delta(p^{-\frac{j}{2}-1}\sfrac{z_2}{z_1})W^{j+2}(p^{-1}z_2)\Bigr)\n
  \!\!\!\!\!\!\!\!\!\!&&
  -\frac{(1-q)(1-t^{-1})}{1-p}
  \Bigl(\delta(p^{\frac{j}{2}}\sfrac{z_2}{z_1})
  \no W^1(p^{-\frac{1}{2}}z_1)W^{j+1}(p^{\frac{1}{2}}z_2)\no\n
  \!\!\!\!\!\!\!\!\!\!&&\qquad\qquad\qquad\qquad
  -\delta(p^{-\frac{j}{2}}\sfrac{z_2}{z_1})
  \no W^1(p^{\frac{1}{2}}z_1)W^{j+1}(p^{-\frac{1}{2}}z_2)\no \Bigr)\n
  \!\!\!\!\!\!\!\!\!\!&&
  +\frac{(1-q)^2(1-t^{-1})^2}{(1-p)^2}
  \biggl(\delta(p^{\frac{j}{2}}\sfrac{z_2}{z_1})
  \Bigl(\frac{p^2}{1-p^2}W^{j+2}(pz_2)+\frac{1}{1-p^j}W^{j+2}(z_2)\Bigr)
  \label{W2Wj}\\
  \!\!\!\!\!\!\!\!\!\!&&\hspace{32mm}
  -\delta(p^{-\frac{j}{2}}\sfrac{z_2}{z_1})
  \Bigl(\frac{p^j}{1-p^j}W^{j+2}(z_2)+\frac{1}{1-p^2}W^{j+2}(p^{-1}z_2)\Bigr)
  \biggr).\nonumber
\ea
Here a normal ordering for currents $\no *\no$ is defined by
\ba
  \!\!\!\!\!\!\!\!\!\!&&\!\!
  \no W^i(rz)W^j(z)\no\n
  \!\!\!\!\!\!\!\!\!\!&=&
  \oint\frac{dz'}{2\pi iz'}\biggl(
  \frac{1}{1-\frac{rz}{z'}}f^{i,j}(\sfrac{z}{z'})W^i(z')W^j(z)
  +\frac{\frac{z'}{rz}}{1-\frac{z'}{rz}}W^j(z)W^i(z')f^{j,i}(\sfrac{z'}{z})
  \biggr)\n
  \!\!\!\!\!\!\!\!\!\!&=&
  \sum_{n\in\Z}\sum_{m=0}^{\infty}\sum_{\ell=0}^m f^{i,j}_\ell \Bigl(
  r^{m-\ell}W^i_{-m}W^j_{n+m}+r^{\ell-m-1}W^j_{n-m-1}W^i_{m+1}
  \Bigr)\cdot z^{-n},
  \label{no}
\ea
where $\frac{1}{1-z}$ stands for $\sum_{n\geq 0}z^n$.
Due to this normal ordering, infinite sums in the RHS of \eq{fWW} become
finite sums on the highest weight representation space.

Eqs.\eq{W1Wj} and \eq{W2Wj} are directly calculated by using the 
commutation relation of $h^i_n$. In principle, we can continue this 
calculation for $i\geq 3$ cases, but in practice it is hopeless. 
So we use another method: fusion and induction.
To write down general formula, we extend the range($0\leq i\leq N$) of 
$W^i(z)$ and that($1\leq i,j\leq N-1$) of $f^{i,j}(z)$ to 
$i\in\Z$ and $i,j\in\Z$ respectively ;
$W^i(z)=0$ for $i<0$ or $i>N$, and $f^{i,j}(z)$ is given by \eq{fij}
for all $i,j\in\Z$.

Explicit expression of the defining relation of $\DWA_{q,t}$ is as follows:
\ba
  \!\!\!\!\!&&f^{i,j}(\sfrac{z_2}{z_1})W^i(z_1)W^j(z_2)
  -W^j(z_2)W^i(z_1)f^{j,i}(\sfrac{z_1}{z_2})\qquad(0\leq i\leq j\leq N)\n
  \!\!\!\!\!&=&-\frac{(1-q)(1-t^{-1})}{1-p}\sum_{k=1}^i\prod_{l=1}^{k-1}
  \gamma(p^{l+\frac12})\n
  \!\!\!\!\!&&\qquad\times\Bigl(\delta(p^{\frac{j-i}{2}+k}\sfrac{z_2}{z_1})
  f^{i-k,j+k}(p^{-\frac{j-i}{2}})W^{i-k}(p^{-\frac{k}{2}}z_1)
  W^{j+k}(p^{\frac{k}{2}}z_2)\n
  \!\!\!\!\!&&\qquad\quad-\delta(p^{-(\frac{j-i}{2}+k)}\sfrac{z_2}{z_1})
  f^{i-k,j+k}(p^{\frac{j-i}{2}})W^{i-k}(p^{\frac{k}{2}}z_1)
  W^{j+k}(p^{-\frac{k}{2}}z_2)\Bigr),
  \label{WiWj}
\ea
where $\ds\gamma(p^{\frac12}z)=\frac{(1-qz)(1-t^{-1}z)}{(1-z)(1-pz)}$.
We can rewrite the RHS of this relation in terms of the normal ordering 
$\no *\no$ by repeated use of the following formula,
which is obtained from \eq{no} and \eq{WiWj},
\ba
  &&f^{i,j}(r^{-1})W^i(rz)W^j(z)\qquad(0\leq i\leq j\leq N)\n
  &=&\no W^i(rz)W^j(z)\no
%
  +\frac{(1-q)(1-t^{-1})}{1-p}
  \sum_{k=1}^i\prod_{l=1}^{k-1}\gamma(p^{l+\frac12})\n
  &&\qquad\quad\times\biggl(
  \frac{1}{1-rp^{-(\frac{j-i}{2}+k)}}f^{i-k,j+k}(p^{-\frac{j-i}{2}})
  W^{i-k}(p^{\frac{j-i+k}{2}}z)W^{j+k}(p^{\frac{k}{2}}z)\n
  &&\qquad\qquad
  -\frac{1}{1-rp^{\frac{j-i}{2}+k}}f^{i-k,j+k}(p^{\frac{j-i}{2}})
  W^{i-k}(p^{-\frac{j-i+k}{2}}z)W^{j+k}(p^{-\frac{k}{2}}z)\biggr),
  \label{noWW}
\ea
where $r\in\mathbb{C}$ is a ``good" number (such that it does not give poles, 
see \eq{fWWpole}).
For example, \eq{W2Wj} is easily recovered by \eq{WiWj} with $i=2$ and 
\eq{noWW} with $i=1$.

In order to prove \eq{WiWj} we present some formulas. Direct calculation shows
\ba
  &&f^{1,j}(p^{\pm\frac{i+1}{2}}z)f^{i,j}(z)=f^{i+1,j}(p^{\pm\frac12}z)
  \times\biggl\{
  \begin{array}{lll}
  1&&i<j\\ \gamma(p^{\pm\frac{i-j+1}{2}}z)&&i\geq j
  \end{array}
  \quad(j\geq 1),
  \label{ff=f}\\
  &&f^{1,i}(p^{\pm(\frac{j-i}{2}+k)}z)f^{1,j}(z)=
  f^{1,i-k}(p^{\pm\frac{j-i+k}{2}}z)f^{1,j+k}(p^{\pm\frac{k}{2}}z)\n
  &&\hspace{73mm}(i,j,i-k,j+k\geq 1),\\
  &&f^{1,i}(p^{\pm\frac{j+i}{2}}z)f^{1,j}(z)=
  f^{1,j+i}(p^{\pm\frac{i}{2}}z)\gamma(p^{\pm\frac{j}{2}}z)
  \quad(i,j\geq 1),
\ea
and $f^{a,b}$ in the RHS of \eq{WiWj} is regular.
By computing $\langle{\lambda}|f^{i,j}(\sfrac{z_2}{z_1})W^i(z_1)W^j(z_2)
\ket{\lambda}$ in the free field realization, we can show that \eq{WiWj}
implies
\ba
  &&\mbox{Poles of $f^{i,j}(\sfrac{z_2}{z_1})W^i(z_1)W^j(z_2)$ 
  ($0\leq i\leq j\leq N$) are}\n
  &&\mbox{$\sfrac{z_2}{z_1}=p^{\pm(\frac{j-i}{2}+k)}$ 
  ($1\leq k\leq\min(i,N-j))$},
  \label{fWWpole}
\ea
because $\langle{\lambda}|f^{i,j}(\sfrac{z_2}{z_1})W^i(z_1)W^j(z_2)
\ket{\lambda}$ is a Taylor series in $\frac{z_2}{z_1}$, and for any states 
of the highest weight representation space, $\ket{\psi}$ and $\ket{\phi}$, 
$\langle\psi|f^{i,j}(\sfrac{z_2}{z_1})W^i(z_1)W^j(z_2)\ket{\phi}$ differs from
$\langle{\lambda}|f^{i,j}(\sfrac{z_2}{z_1})W^i(z_1)W^j(z_2)\ket{\lambda}$
only for finite number of terms(Laurant polynomials in $z_1$ and $z_2$), 
which do not create other poles.
(See also Appendix C of ref.\cite{HJKOS99} where different notation is used.)
Therefore $f^{a,b}W^aW^b$ in the RHS of \eq{WiWj} is regular and we can 
reverse its order, $f^{a,b}(p^c)W^aW^b=f^{b,a}(p^{-c})W^bW^a$.
{}From \eq{W1Wj} (or by using the free field realization), 
we have the following fusion relation
\ba
  &&\lim_{z_1\rightarrow p^{\pm\frac{j+1}{2}}z_2}\bigl(1-p^{\pm\frac{j+1}{2}}
  \sfrac{z_2}{z_1}\bigr)f^{1,j}(\sfrac{z_2}{z_1})W^1(z_1)W^j(z_2)\n
  &=&\mp\frac{(1-q)(1-t^{-1})}{1-p}W^{j+1}(p^{\pm\frac12}z_2)
  \quad(1\leq j\leq N),
  \label{fusionW1Wj}
\ea
and if \eq{WiWj} is correct, \eq{WiWj} implies
\ba
  &&\lim_{z_2\rightarrow p^{\mp\frac{j+i}{2}}z_1}\bigl(1-p^{\mp\frac{j+i}{2}}
  \sfrac{z_1}{z_2}\bigr)f^{j,i}(\sfrac{z_1}{z_2})W^j(z_2)W^i(z_1)\n
  &=&\pm\frac{(1-q)(1-t^{-1})}{1-p}\prod_{l=1}^{i-1}\gamma(p^{l+\frac12})\cdot
  W^{j+i}(p^{\mp\frac{j}{2}}z_1)\quad(0\leq i\leq j\leq N).
  \label{fusionWiWj}
\ea

\noindent
{\bf Proof of \eq{WiWj} : }
(i) The case $i=0$ and $i\leq\forall j\leq N$, 
and the case $j=N$ and $0\leq\forall i\leq j$ are trivial.  
(ii) The case $i=1$ and $i\leq\forall j\leq N$, i.e. \eq{W1Wj}, is already 
proved. 
(iii) Let us assume \eq{WiWj} holds for $i(<N)$ and $i\leq\forall j\leq N$.
We will show \eq{WiWj} holds for $i+1$ and $i+1\leq\forall j\leq N$.
(For $i=N-1$, we have $i+1=N\leq j=N$. Therefore it is sufficient to 
consider $i<N-1$ and $j<N$.)
Multiply $f^{1,i}(\sfrac{z_1}{z_3})f^{1,j}(\sfrac{z_2}{z_3})W^1(z_3)$ from
left to \eq{WiWj} with $i\geq 1$ (whose second $f^{a,b}(p^c)W^aW^b$ term 
in the RHS is replaced by reversed order one $f^{b,a}(p^{-c})W^bW^a$),
rewrite $f^{1,j}(\sfrac{z_2}{z_3})W^1(z_3)W^j(z_2)
=f^{j,1}(\sfrac{z_3}{z_2})W^j(z_2)W^1(z_3)+\cdots$ by using \eq{W1Wj},
multiply $\frac{1-p}{(1-q)(1-t^{-1})}(1-p^{-\frac{i+1}{2}}\sfrac{z_2}{z_1})$,
and take a limit $z_3\rightarrow p^{-\frac{i+1}{2}}z_1$.
By using \eq{ff=f}--\eq{fusionW1Wj} and \eq{fusionWiWj}
(with $j\rightarrow j+1$), studying poles carefully and replacing 
$z_1\rightarrow p^{\frac12}z_1$, we obtain \eq{WiWj} with $i\rightarrow i+1$ 
($2\leq i+1\leq j<N$). 
(\romannumeral4) Therefore we have proved \eq{WiWj} by induction on $i$.
\qed

\section{Relation to $\cal{Z}$-Algebra}

Affine Lie algebra $\widehat{\mathfrak{sl}}_N$ is an associative algebra
over $\mathbb{C}$ with the Chevally generators, $e^{\pm}_i$ and $h_i$ 
($i=0,1,\cdots,N-1$), which satisfy
\be
  [h_i,h_j]=0,\quad[h_i,e^{\pm}_j]=\pm a_{ij}e^{\pm}_j,\quad
  [e^+_i,e^-_j]=\delta_{ij}h_i,
\ee
and the Serre relation $\ds({\rm ad}\,e^{\pm}_i)^{1-a_{ij}}e^{\pm}_j=0$ 
($i\neq j$), where $(a_{ij})_{0\leq i,j\leq N-1}$ is the Cartan matrix 
of $A^{(1)}_{N-1}$ Dynkin diagram. \cite{K90}
This algebra admits various gradations and we denote its grading operator 
as $d$ and $\rho$ for the homogeneous and principal gradation respectively, 
which satisfy
\be
  \begin{array}{lll}
  \mbox{homogeneous gradation}&:&[d,e^{\pm}_i]=\pm e^{\pm}_i\delta_{i0},\\
  \mbox{principal gradation}&:&[\rho,e^{\pm}_i]=\pm e^{\pm}_i.
  \end{array}
\ee 
In current basis $\widehat{\mathfrak{sl}}_N$ is given as follows:\\
{\bf homogeneous gradation}\\
\underline{generators} : $H^i_n$, $E^{\pm,i}_n$ ($n\in\Z$, $1\leq i\leq N-1$),
$k$ : center, $d$ : grading operator.\\ 
\underline{relations} :
\ba
  &&[H^i_n,H^j_m]=k\bar{a}_{ij}n\delta_{n+m,0},\quad
  [H^i_n,E^{\pm,j}_m]=\pm\bar{a}_{ij}E^{\pm,j}_{n+m},\n
  &&[E^{+,i}_n,E^{-,j}_m]=\delta^{ij}(H^i_{n+m}+kn\delta_{n+m,0}),\quad
  [d,X_n]=nX_n\;\;(X=H^i,E^{\pm,i})\,,
\ea
and $[E^{\pm,i}_n,E^{\pm,j}_m]=[E^{\pm,i}_{n-1},E^{\pm,j}_{m+1}]$ and
the Serre relations which we omit to write explicitly, 
where $(\bar{a}_{ij})_{1\leq i,j\leq N-1}$ is the Cartan matrix 
of $A_{N-1}$ Dynkin diagram. 

\noindent
{\bf principal gradation} $\quad$ Let us set $\omega=e^{\frac{2\pi i}{N}}$.
Symbol $\equiv$ stands for $\equiv\!\!\!\!\!\!\pmod{N}$. \\
\underline{generators} : $\beta_n$ ($n\in\Z$, $n\not\equiv 0$), 
$x^{(\mu)}_n$ ($n\in\Z$, $1\leq\mu\leq N-1$, $\mu$ is understood as\\
\phantom{\underline{generators} : } 
$\bmod$ $N$), $k$ : center, $\rho$ : grading operator.\\ 
\underline{relations} :
\ba
  &&\lbrack\beta_n,\beta_m\rbrack=kn\delta_{n+m,0}\quad(n,m\not\equiv 0)\,,
  \quad
  \lbrack\beta_n,x^{(\nu)}_m\rbrack=(1-\omega^{-\nu n})x^{(\nu)}_{n+m}\,,\n
  &&\lbrack x^{(\mu)}_n,x^{(\nu)}_m\rbrack=
  \biggr\{\begin{array}{ll}
  (\omega^{-\mu m}-\omega^{-\nu n})x^{(\mu+\nu)}_{n+m}&(\mu+\nu\not\equiv 0)\\
  (\omega^{-\mu m}-\omega^{-\nu n})\beta_{n+m}+kn\omega^{\mu n}\delta_{n+m,0}
  &(\mu+\nu\equiv 0)\,,
  \end{array}
  \label{principalslN}\\
  &&[\rho,X_n]=X_n\;\;(X=\beta,x^{(\mu)})\,,\nonumber
\ea
and the Serre relations.\\
Since these two current basis are basis of the same Lie algebra 
$\widehat{\mathfrak{sl}}_N$, they are related by linear transformation,
\ba
  &&\beta_{Nm+\nu}=\sum_{i=1}^{N-\nu}E^{i,i+\nu}_m
  +\sum_{i=N-\nu+1}^NE^{i,i+\nu-N}_{m+1},\n
%
  &&x^{(\mu)}_{Nm+\nu}=\sum_{i=1}^{N-\nu}\omega^{\mu(i+\nu-1)}E^{i,i+\nu}_m
  +\sum_{i=N-\nu+1}^N\omega^{\mu(i+\nu-1)}E^{i,i+\nu-N}_{m+1},\\
  &&x^{(\mu)}_{Nm}=\sum_{i=1}^{N-1}\frac{1-\omega^{\mu i}}{1-\omega^{\mu}}H^i_m
  -\frac{k}{1-\omega^{\mu}}\delta_{m,0},\nonumber
\ea
where $m\in\Z$ and $1\leq\mu,\nu\leq N-1$. Here, for simplicity of the
presentation, we have introduced $\widehat{\mathfrak{gl}}_N$ generators 
$E^{i,j}_n$ ($n\in\Z$, $1\leq i,j\leq N$), which satisfy 
$[E^{i,j}_n,E^{i',j'}_m]=\delta^{ji'}E^{i,j'}_{n+m}
-\delta^{ij'}E^{i',j}_{n+m}+\delta^{ij'}\delta^{ji'}kn\delta_{n+m,0}$,
and the generators in the homogeneous picture are expressed as 
$E^{+,i}_n=E^{i,i+1}_n$, $E^{-,i}_n=E^{i+1,i}_n$ 
and $H^i_n=E^{i,i}_n-E^{i+1,i+1}_n$.
We remark $E^{i,j}_n=[E^{i,l}_m,E^{l,j}_{n-m}]$ ($i\neq j$) and
this RHS is independent on $l$ and $m$.

Next let us consider the splitting of the Cartan part:
\be
  (\mbox{ $\widehat{\mathfrak{sl}}_N$ generator})=
  (\mbox{exponential of Cartan generators})\times
  (\mbox{new generator}),
\ee
where new generator commutes with Cartan generators.
For homogeneous gradation, the algebra generated by these new generators 
is known as the ($\mathfrak{sl}_N$ version of) parafermion algebra of 
level $k$. For principal gradation, we name it as ($\mathfrak{sl}_N$-version
of) ${\cal Z}$ algebra of level $k$, $\ZA$. ($N=2$ case was studied by
Lepowsky and Wilson.\cite{LW81})
Explicitly the generators of $\ZA$, $z^{\mu}_n$ ($n\in\Z$, 
$1\leq\mu\leq N-1$, $\mu$ is understood as $\bmod$ $N$), are obtained by
\be
  x^{(\mu)}(\zeta)=\,:\exp\Bigl(-\frac{1}{k}\sum_{n\not\equiv 0}\frac{1}{n}
  (1-\omega^{\mu n})\beta_n\zeta^{-n}\Bigr): z^{\mu}(\zeta),
\ee
where $:*:$ stands for the normal ordering for boson $\beta_n$ and we have 
introduced currents $x^{(\mu)}(\zeta)=\sum_{n\in\Z}x^{(\mu)}_n\zeta^{-n}$ and
$z^{\mu}(\zeta)=\sum_{n\in\Z}z^{\mu}_n\zeta^{-n}$.
Then \eq{principalslN} implies the relation of $\ZA$,
\ba
  &&g^{\mu,\nu}(\sfrac{\zeta_2}{\zeta_1})z^{\mu}(\zeta_1)z^{\nu}(\zeta_2)
  -z^{\nu}(\zeta_2)z^{\mu}(\zeta_1)g^{\nu,\mu}(\sfrac{\zeta_1}{\zeta_2})\n
  &=&\biggr\{\begin{array}{ll}
  \delta(\omega^{\mu}\sfrac{\zeta_2}{\zeta_1})z^{\mu+\nu}(\omega^{\mu}\zeta_2)
  -\delta(\omega^{-\nu}\sfrac{\zeta_2}{\zeta_1})z^{\mu+\nu}(\zeta_2)
  &(\mu+\nu\not\equiv 0)\\
  k\,D\delta(\omega^{\mu}\sfrac{\zeta_2}{\zeta_1})&(\mu+\nu\equiv 0)\,,
  \end{array}
  \label{ZA}
\ea
where $D=\zeta\frac{d}{d\zeta}$, $D\delta(\zeta)=\sum_{n\in\Z}n\zeta^n$
and the structure function $g^{\mu,\nu}(\zeta)$ is given by
\be
  g^{\mu,\nu}(\zeta)=\exp\Bigl(-\frac{1}{k}\sum_{n>0\atop n\not\equiv 0}
  \frac{1}{n}(1-\omega^{\mu n})(1-\omega^{-\nu n})\zeta^n\Bigr).
\ee

Next we present an interesting relation between $\DWA_{q,t}$ and $\ZA$.
Let us consider the following limit\footnote{
For this choice of $t=\omega^{-1}q^{\frac{k+N}{N}}$, 
we cannot take Limit I because $\beta=\frac{k+N}{N}-\frac{2\pi i}{N\hbar}$ 
depends on $\hbar$.
}:
\be
  \mbox{Limit II : }
  \biggl\{\begin{array}{ll}
  q=e^{\hbar}\,,&\hbar\rightarrow 0\\
  t=\omega^{-1}q^{\frac{k+N}{N}}\,,\quad&k\mbox{ : fixed}\,.
  \end{array}
\ee
We assume that DWA currents $W^i(z)$ have the $\hbar$-expansion
\be
  W^i(p^{\frac{1-i}{2}}\zeta)=\hbar\omega^{\frac{i}{2}}z^i(\zeta)+O(\hbar^2).
  \label{Wz}
\ee
Then, under the Limit II, we can show that the relation of $\DWA_{q,t}$ 
\eq{WiWj} reduces to that of $\ZA$ \eq{ZA}.
(Eq.\eq{WiWj} begins from $\hbar^2$ term and its coefficient is \eq{ZA}.
We remark that in this derivation we do not use free field realization 
at all.) 
In other words, $\DWA_{q,t}$ with $t=\omega^{-1}q^{\frac{k+N}{N}}$ can be
considered as a $q$-deformation of $\ZA$, which we denote as $\DZA$,
\be
  \DZA=\DWA_{q,t=\omega^{-1}q^{\frac{k+N}{N}}}.
\ee

Concerning the free field realization, however, our assumption \eq{Wz} does
not hold on the Fock space except for $N=2$ case. 
But calculation of some correlation functions supports the assumption \eq{Wz};
We have checked $\langle\lambda|W^1(\zeta_1)\cdots W^1(\zeta_n)\ket{\lambda}
=O(\hbar^n)$ for $n\leq 6$. We guess that the assumption \eq{Wz} holds on 
the level of correlation functions, or, on the irreducible representation 
space obtained by taking some BRST cohomology. For $N=2$ case, \eq{Wz} holds
on the Fock space, and screening currents and vertex operators of DVA (after
some modification of zero mode) reduce to those of ZA.

Finally we mention the character of ${\rm DZA}(\widehat{\mathfrak{sl}}_2)_k$ 
for $k\in\Z_{\geq 2}$, i.e. 
that of ${\rm DVA}_{q,t}={\rm DWA}(\widehat{\mathfrak{sl}}_2)_{q,t}$ with 
$t=e^{-\pi i}q^{\frac{k+2}{2}}$. 
We write $W^1(\zeta)$ and $w^1(\lambda)$ as $T(\zeta)$ and $\lambda$ 
respectively, e.g., the highest weight state is defined by 
$T_n\ket{\lambda}=\lambda\ket{\lambda}\delta_{n0}$ ($n\geq 0$). 
Since degenerate representations of DVA occur at 
$\lambda=\lambda_{r,s}=t^{\frac{r}{2}}q^{-\frac{s}{2}}
+t^{-\frac{r}{2}}q^{\frac{s}{2}}$,\cite{O99} let us consider 
$\lambda=\lambda_{1,j+\frac{k+2}{2}}$ ($j=-\frac{k}{2},-\frac{k}{2}+1,\cdots,
\frac{k}{2}$) representations.
Grading operator $\rho$ satisfies $[\rho,T_n]=nT_n$ and $-\rho\ket{\lambda}=
\bigl(\frac{2j^2+k}{4(k+2)}-\frac18\bigr)\ket{\lambda}$.
The character of DZA is defined by $\chi^{{\rm DZA}}_j(\tau)={\rm tr}\,
y^{-\rho}$ where $y=e^{2\pi i\tau}$ and the trace is taken over irreducible 
DZA spin $j$ representation space.
Shiraishi and present author conjectured\footnote{
We remark that this character appears in the calculation of the
one-point local height probability of the Kashiwara-Miwa model (
M.~Jimbo, T.~Miwa and M.~Okado, 
\reftitle{Solvable Lattice Models with Broken $\Z_N$ Symmetry and Hecke's
Indefinite Modular Forms}
\refj{Nucl. Phys.} {\bf B275[FS17]} (1986) 517-545
).
} 
\be
  \chi^{{\rm DZA}}_j(\tau)=y^{\frac{2j^2+k}{4(k+2)}-\frac18}
  \frac{1}{(y;y)_{\infty}}\sum_{m\in\Z}(-1)^my^{m(j+\frac{k+2}{2}m)}
  =y^{\frac{2j^2+k}{4(k+2)}-\frac18}\chi^{(2,k+2)}_{1,j+\frac{k+2}{2}}(\tau).
  \label{char}
\ee
Here $\chi^{(p',p'')}_{r,s}(\tau)$ is the Rocha-Caridi character 
formula\cite{O99}
\be
  \chi^{(p',p'')}_{r,s}(\tau)=\frac{1}{(y;y)_{\infty}}\sum_{m\in\Z}
  \Bigl(y^{(p''r-p's+mp'p'')m}-y^{(r+mp')(s+mp'')}\Bigr),
\ee
which gives the character of the Virasoro minimal representation 
when $p'$ and $p''$ are coprime, $(p',p'')=1$. In the above case,
$p'=2$ and $p''=k+2$ imply that $(p',p'')=1$ for odd $k$ but 
$(p',p'')=2$ for even $k$.
When $q$ is not a root of unity, by studying the Kac determinant of 
DVA\cite{O99}, we can check that \eq{char} is true. 
We remark that the character of ${\rm DZA}(\widehat{\mathfrak{sl}}_2)_k$, 
$\chi^{{\rm DZA}}_j$, coincides with that of 
${\rm ZA}(\widehat{\mathfrak{sl}}_2)_k$ which is obtained by using the result 
of ref.\cite{HJM00}, BRST structure of principal $\widehat{\mathfrak{sl}}_2$.

\section*{Acknowledgements}

The author would like to thank 
H.~Awata and J.~Shiraishi
for valuable discussions,
and M.~Jimbo for helpful comments
and explaining poles of $f^{i,j}W^iW^j$ and efficient graphical computational
method for DWA correlation functions. 
He is also grateful to M.L.~Ge for kind invitation to this APCTP-Nankai 
joint symposium and thanks members of Nankai Institute for their kind 
hospitality. 
This work is supported in part by Grant-in-Aid for Scientific Research from 
the Ministry of Education, Culture, Sports, Science and Technology,
No.12640261.

\appendix
\section{$\DWA_{q,t}$ and $\zeta$-function regularization\\
(by H.~Awata, S.~Odake and J.~Shiraishi)}

In this appendix we present an interesting relation between $\DWA_{q,t}$ and 
$\zeta$-function regularization.\cite{AOS00}

In string theory\cite{string}, the physical state condition is given by
$(L_0-1)\ket{{\rm phys}}=0$ (and its antichiral counterpart),
where $L_0$ is the zero mode of the Virasoro generator.
This condition and the space-time dimension are derived by careful
study of string theory (Lotentz invariance in the light-cone gauge,
nilpotency of BRST charge, etc.), but there is a shortcut method, 
$\zeta$ function regularization method.

First we illustrate this method by taking a bosonic string theory 
as an example. In the light-cone gauge the Virasoro generator $L_n$ is
given by 
$L_n=\sum_{i=1}^{24}\sum_{m\in\Z}\sfrac{1}{2}:\alpha^i_{n-m}\alpha^i_m:$
where $\alpha^i_n$ ($n\in\Z$, $i=1,\cdots,24$) satisfies 
$[\alpha^i_n,\alpha^j_m]=n\delta^{ij}\delta_{n+m,0}$ and $:*:$ stands for
the normal ordering.
The Virasoro zero mode without the normal ordering is 
$L_0^{{\rm noNO}}
=\sum_{i=1}^{24}\sum_{n\in\Z}\sfrac{1}{2}\alpha^i_{-n}\alpha^i_n
=\sum_{i=1}^{24}\sum_{n\in\Z}\sfrac{1}{2}:\alpha^i_{-n}\alpha^i_n:
+12``\sum_{n>0}n"$.
Of course the sum $``\sum_{n>0}n"$ is divergent and this expression is
meaningless. But we replace the sum $``\sum_{n>0}n"$ by $\zeta(-1)$, 
where $\zeta(z)$ is the Riemann $\zeta$ function. 
Then the above physical state condition is equivalent to the condition that
the Virasoro zero mode without the normal ordering annihilates the physical
state:
\be
  L_0^{{\rm noNO}}\ket{{\rm phys}}=0,\quad
  L_0^{{\rm noNO}}=L_0+12\zeta(-1),
\ee
because of the value $\zeta(-1)=-\frac{1}{12}$.
We might say that the Virasoro generator ``knows" the value $\zeta(-1)$.
 
Next let us mimic the above procedure for $\DWA_{q,t}$ case.
DWA current without the normal ordering becomes
$W^{i\,\,{\rm noNO}}(z)=``f^{i,i}(1)^{-\frac12}"\,W^i(z)$, 
where $``f^{i,i}(1)"$ is divergent for generic $\beta$ (recall $t=q^{\beta}$ 
and $q=e^{\hbar}$). 
Let $a^i_{2m}$ be coefficients of the following $\hbar$-expansion 
$(1-q^n)(1-t^{-n})\frac{1-p^{in}}{1-p^n}\frac{1-p^{(N-i)n}}{1-p^{Nn}}
=\sum_{m>0}a^i_{2m}(n\hbar)^{2m}$. 
Then $f^{i,i}(z)$ is $f^{i,i}(z)=\exp\bigl(\sum_{n>0}\frac{1}{n}
\sum_{m>0}a^i_{2m}(n\hbar)^{2m}z^n\bigr)$.
We define $\zeta$-regularized $f^{i,i}_{\zeta\mbox{-}{\rm reg}}(1)$ by 
exchanging these summations over $n$ and $m$ and 
replacing $\sum_{n>0}n^{2m-1}$ with $\zeta(1-2m)$ as follows:
\be
  f^{i,i}_{\zeta\mbox{-}{\rm reg}}(1)=\exp\Bigl(\sum_{m>0}a^i_{2m}\zeta(1-2m)
  \hbar^{2m}\Bigr).
\ee
In the Limit I \eq{LimitI}, DWA current behaves
as $W^i(z)={N\choose i}+O(\hbar^2)$, which can be shown by using free
field realization.\cite{AKOS95} So we require that $\beta=\frac{N+1}{N}$ or 
$\frac{N}{N+1}$, which corresponds to the vanishing Virasoro central charge,
and the zero mode of the $i$-th DWA current without normal ordering takes
the above value ${N\choose i}$ on the vacuum state $\ket{{\rm vac}}$, which
is characterized by $h^i_n\ket{{\rm vac}}=0$ ($n\geq 0$, $\forall i$),
\be
  W^{i\,\,{\rm noNO}}_0\ket{{\rm vac}}
  =\mbox{$\Bigl({N\atop i}\Bigr)$}\ket{{\rm vac}},\quad
  W^{i\,\,{\rm noNO}}(z)
  =f^{i,i}_{\zeta\mbox{-}{\rm reg}}(1)^{-\frac12}\,W^i(z).
\ee
Since we can show $W^i_0\ket{{\rm vac}}=\Bigl[{N\atop i}\Bigr]\ket{{\rm vac}}$,
this requirement implies
\be
  f^{i,i}_{\zeta\mbox{-}{\rm reg}}(1)^{\frac12}
  =\mbox{$\Bigl({N\atop i}\Bigr)^{-1}$}\mbox{$\Bigl[{N\atop i}\Bigr]$},
\ee
where $\Bigl[{N\atop i}\Bigr]=\frac{[N]!}{[i]![N-i]!}$, $[n]!=[n]\cdots[1]$ 
and $[n]=\frac{p^{\frac{n}{2}}-p^{-\frac{n}{2}}}{p^{\frac12}-p^{-\frac12}}$.
We can check that this equation really holds by using formulas
$\log(\sinh x)=\log x+\sum_{n>0}(-1)^{n-1}\frac{2^{2n-1}B_n}{(2n)!n}$ $x^{2n}$ 
($0<|x|<\pi$) and $\zeta(1-2m)=(-1)^m\frac{B_m}{2m}$ ($m=1,2,\cdots$).
Here $B_n$ is the Bernoulli number defined by $\frac{x}{e^x-1}+\frac{x}{2}=
1+\sum_{n>0}(-1)^{n-1}\frac{B_n}{(2n)!}x^{2n}$ ($|x|<2\pi$).
Therefore we might say that $\DWA_{q,t}$ (with 
$t=q^{\frac{N+1}{N}},q^{\frac{N}{N+1}}$) for each $N$ ``knows" all the values 
$\zeta(1-2m)$ ($m=1,2,\cdots$).


\end{document}